\documentclass{amsart}
\usepackage[english]{babel}
\usepackage{amsmath,amsthm,amssymb,verbatim}
\usepackage{hyperref}
\usepackage{graphicx}
\usepackage{color}

\newtheorem{theorem}{Theorem}

\newtheorem{definition}{Definition}

\newtheorem{statement}{Proposition}

\newtheorem{problem}{Problem}
\newtheorem{definition-theorem}{Theorem-Definition}

\begin{document}

\title{Classification of finite metric spaces and combinatorics of convex polytopes}

\author{A.~M.~Vershik}

\thanks{Supported by the Russian Science Foundation grant 14-11-00581.}

\address{St.~Petersburg Department of Steklov Institute of Mathematics and St.~Petersburg State University}

\email{vershik@pdmi.ras.ru}

\maketitle
\rightline{\emph{To the memory of Dima Arnold}}

\begin{abstract}

We describe the canonical correspondence between finite metric spaces and symmetric convex polytopes, and formulate the problem about classification of the metric spaces in terms of combinatorial structure of those polytopes.

 \end{abstract}

  \section{Introductory remark}

  We will discuss problems with very elementary formulations that concern
  the most popular notions in mathematics: metric spaces, convex geometry,
  combinatorics of polytopes and Kantorovich's optimal transportation.
  According to Arnold's classification, there are two ways to introduce a new subject:
  the first way (he called it the ``Russian tradition'') is to start with
  ``the simplest non-trivial partial case''. ---
  I will use this approach.
  The second and the opposite tradition, which I also like very much
  (he called it ``Bourbaki's tradition'') is to start with an
  ``extremely general case that is impossible to further generalize''.

 So my metric spaces will be finite, polytopes will be finite-dimensional etc. but all the notions and problems
 have infinite, infinite-dimensional, and continuous analogs.

 \textbf{GENERAL PROBLEM.}
 Study and classify finite metric spaces according to combinatorial
 properties of their fundamental polytopes, associated with metric spaces in
 a canonical fashion.
 A more precise formulation follows.

  In the rest of the paper, I explain the terminology that is used in the title and in the Problem.

  \section{Finite metric spaces, canonical polytopes, and problems}

1. Let $(X,\rho)$, $|X|=n$, be a finite metric space.
We will write $V(X)$ for the vector space of all real-valued functions on $X$.
The value of a function $v\in V(X)$ at a point $x\in X$ will be denoted by $v^x$.
The space $V(X)$ can be naturally identified with the space of all formal
linear combinations of elements of $X$.
Under this identification, any element $x\in X$ identifies with
the delta-function equal to 1 at $x$ and to 0 at all other points of $X$.
Define the real vector space $V_0(X)$ as a subspace of $V(X)$
consisting of all $v\in V(X)$ with $\sum_{x\in X}v^x=0$.

Let us define the map $\delta:X\to V_0(X)$ taking an element $x$
to the function $\delta(x)=\delta_x\in V_0(X)$ such that $\delta_x^y=-\frac 1n$
for all $y\ne x$.
Then we must have $\delta_x^x=\frac{n-1}n$.

The convex hull $\mathrm{Conv}[\{\delta(x)\}, x\in X]$
of the image of this map is a $(n-1)$-dimensional simplex
denoted by $\Sigma(X) \subset V_0(X)$
(this simplex is obtained from the standard coordinate simplex in $V(X)$
by the projection mapping $x\in X$ to $\delta(x)$).

The metric $\rho$ can be considered now as a metric on the vertices
of the simplex $\Sigma(X)$,
we use the same symbol $\rho$ to denote this metric.
For every pair of distinct points $x$, $y\in X$, consider the vectors
$e_{x,y} \in V_0$ defined by the formula:
$$
e_{x,y}=\frac{\delta_x-\delta_y}{\rho(x,y)}.
$$
This vectors will play a major role in what follows.

 \begin{definition}
 \textbf{The fundamental polytope of a metric space} $(X,\rho)$ is the convex polytope
 $R_{X,\rho}$ obtained as the convex hull of all vectors $e_{x,y}$, where $(x,y)$
 run through all pairs of distinct points of $X$.
 The combinatorial type of $R_{X,\rho}$, i.e., the isomorphism class of the
 corresponding poset of faces,
 is called the \emph{combinatorial structure} of the finite metric space $(X,\rho)$.
 Two finite metric spaces with the same combinatorial structures called \emph{similar metric spaces}.
 \end{definition}

It is easy to see that the fundamental polytope $R_{X,\rho}$
is centrally symmetric (i.e., it coincides with its reflection in the origin).
If $\rho(x,y)=1$ for any pair of distinct elements $x$, $y\in X$, then
the fundamental polytope is the Minkowski sum of two simplices $\Sigma(X)$ and
$-\Sigma(X)$.

In general,
we consider the fixed set of rays $\{\lambda (\delta_x -\delta_y): \lambda > 0\}$,
which is independent of the metric, and choose one point
in each ray with $\lambda=\frac{1}{\rho(x,y)}$ (this choice now depends on the
metric).
The fundamental polytope is the convex hull of all thus obtained points.

A simple fact \cite{MPV} that characterizes the fundamental polytopes
of finite metric spaces, is the following:
given a symmetric function $\rho:X\times X\to \mathbb{R}_{>0}$,
consider the polytope
$$
R_{X,\rho}=\mathrm{Conv}\left\{e_{x,y}=\frac{\delta_x-\delta_y}{\rho(x,y)}:
\, x\ne y,\ x,y \in X\right\}.
$$
Then this function $\rho$ is a metric on $X$ if and only if
no point $e_{x,y}$ belongs to the interior of convex hull of the other points.

The polytope $R_{X,\rho}$ is convex and central symmetric; therefore, it
defines the Minkowski-Banach norm $\|\cdot\|_{\rho}$
in the real vector space $V_0(X)$, whose unit ball is by definition
the polytope $R_{X,\rho}$.
In the finite-dimensional case, this norm is the so called
\emph{Kantorovich--Rubinstein norm}.
If $\rho(x,y)=1$ whenever $x\ne y$, then the corresponding fundamental polytope
is the so called \emph{root polytope}, and the corresponding Kantorovich--Rubinstein
norm is the restriction to $V_0(X)$ of the $\ell^1$-norm in the space $V(X)$.

Thus we reduce the analysis of metric space to the convex geometry of fundamental polytope.
Since the space $X$
is isometrically embedded into $V_0(X)$ (points of $X$ correspond to
the vertices of the simplex $\Sigma(X)$, endowed with metric $\rho$ (see above)
we can restore metric on $X$ using fundamental polytope.

Recall that the combinatorial type of a convex polytope is
the isomorphism class of the partial ordered set form by the faces
of the polytope, ordered by inclusion; the $f$-vector of a convex polytope
is the finite tuple $(f_0,f_1, \dots f_n)$, where
$f_0=n$ is the number of vertices, $f_1$ is the number of edges, etc.,
$f_{n-1}$ is number of facets (i.e., faces of codimension 1) and, finally, $f_n=1$.

Our definition is functorial in the sense that each isometry of one metric space to another vector space produces an affine isomorphism of the corresponding fundamental polytopes.

We may say that the notion of the combinatorial type of metric
spaces defines a natural stratification of the cone ${\mathcal M}_n$ of
all distance matrices (i.e., real symmetric $n\times n$ matrices, whose
entries are the values of a distance function defined on some finite metric
space of cardinality $n$).
Below we suggest to study this important stratification, more precisely,
to solve the following problems.

\begin{problem}
1. Express the combinatorial structure of $(X,\rho)$ in terms of the
metric $\rho$, i.e., find the $f$-vector of the corresponding fundamental polytope
in terms of the metric $\rho$ itself --- using linear inequalities on
the values of metric etc.

2. Estimate the number of combinatorial structures for any given $n$ and study its
asymptotic behavior as $n$ tends to infinity.
The most interesting thing is to estimate the number of ``open''
(generic) types.

3. Provide sufficient conditions on two metric spaces to be similar.

It is well known that most finite metric spaces cannot be isometrically
imbedded into a Euclidean or a Hilbert space (we say that these metrics
do not have Euclidean type).
The following question appears:

4. Describe the combinatorial types of metric spaces of Euclidean type.
Do we obtain all combinatorial types?

5. Does the stratification of the space of distance matrices into the
combinatorial types is \emph{universal}?
Or, on the contrary, there are some restrictions on the topological types
of the open components?\footnote{A classification problem (in algebraic geometry, combinatorics, etc.) of a certain set of objects up to a certain equivalence relation is called a ``universal problem'' if all possible kinds of singularities or stratifications, say, of arbitrary semi-algebraic varieties can occur in the topology of equivalence classes.

A well-known example is Mnev's theorem on the universality of the combinatorial
classification of convex polytopes in dimensions $\geq 4$.
(see the papers by A.Vershk and N.Mnev in  ``Topology and geometry --- Rokhlin Seminar''. Lecture Notes in Math. 1346, 557-581 (1988) and more recent literature.}
\end{problem}

\section{The simplest example: a metric on a Cartan subalgebra.}

Consider a very degenerate metric space with $n$ points:
$X=\textbf{n}$, the set of all integers from 1 to $n$ with mutual distances
between all pairs of distinct points equal to $1$:
$\rho(i,j)=\delta_{i,j}$.

In this case, the simplex $\Sigma(X)$ is a regular Euclidean simplex.
We can view it as $(n-1)$-simplex in a Cartan subalgebra of the Lie algebra $A_n$.
From this viewpoint, the simplex is spanned by all positive simple roots
$e_{i,i+1}$, $i=1$, $\dots$, $n-1$, and one maximal root $e_{n,1}$.

\begin{statement}
Let $X=\textbf{n}$ and $\rho(i,j)=\delta_{i,j}$, $i,j=1$, $\dots$, $n$.
Identify the vector space $V_0(X)$ with a Cartan subalgebra of the Lie algebra $A_n$.
Then the fundamental polytope $R_{(X,\rho)}$
is the convex hull of all roots.
The norm $\|.\|_{\rho}$ associated with the fundamental polytope coincides with
the restriction of the $\ell^1$-norm on $V(X)$
$$
\left|v=(v_1,v_2 \dots v_n)\right|_{\ell^1}=\sum_{i=1}^n |v_i|,
$$
to the hyperplane $V_0(X)\subset V(X)$.
Thus, the polytope $R_{X,\rho}$ in this case is the intersection of
an $n$-dimensional octahedron with the hyperplane $V_0(X)$.

The corresponding norm $\|\cdot\|$ on the Lie algebra of
skew-hermitian matrices with zero trace
is the ``nuclear norm'', for which the norm of a matrix is
the sum of the moduli of all its eigenvalues.
\end{statement}

It is natural to call $R_{X,\rho}$ the ``root polytope''.
An easy exercise is to find the $f$-vector in this case.
For example, if $n=3$ and $\dim(V_0)=2$, then the fundamental polygon is a regular hexagon,
and the norm is the hexagonal norm.
For $n=4$, see the figure: the $f$-vector is equal to $(12,24,14)$,
the facets are 8 regular triangles and 6 squares.
For $n=3$, the group of symmetries of the fundamental polytope is the dihedral group
$D_6={\mathbb{Z}}_2\rightthreetimes {\mathbb{Z}}_6$.
For $n=4$, see the figure.

Note that the group of symmetries of the root polytope includes
the Weyl group (which is isomorphic to the symmetric group).
Root polytopes were mentioned for other reason in \cite{GK}.

\section{Why such combinatorial types?
The Kantorovich metric and geometry of the transportation problem}

Why is the definition of the fundamental polytope $R_{X,\rho}$ associated with a
finite metric space $(X,\rho)$ natural?
The justification is as follows.
We want to define a natural metric $\bar \rho$ on $V_0(X)$ that extends the metric
$\rho$ on $X$.
The latter is identified with the set of vertices of the simplex
$\Sigma(X)$ so that the distance between two points $x$ and $y$ in $X$
coincides with the distance between the corresponding vertices $v_x$ and $v_y$ of the simplex $\Sigma(X)$.
In another words, we want to define a norm $\|\cdot\|$ in $V_0(X)$
with the property $\|e_{x,y}\|=\rho(x,y)$.
There are many such extensions but there is a maximal one:

\begin{theorem}(\cite{MPV})
The norm $\|\cdot\|_{\rho}$, called the Kantorovich-Rubinstein norm,
is the unique maximal extension of $\rho$: all other norms possessing
this extension property are less that this one because
the fundamental polytope is contained in the unit balls of all such norms.
\end{theorem}

The direct description of an extension of the metric $\rho$
to the whole simplex $\Sigma(X)$ is a classical definition
by Kantorovich of his transportation metric.
This definition is well known in the mathematical economics
and in linear programming but there were no publications
(before \cite{MPV}, see comments below), in which
characteristic properties of fundamental polytopes are discussed
and serious combinatorial investigations are conducted.
Recall, for the sake of completeness, the classical definition of the
\emph{Kantorovich transportation metric},
which is equivalent to our definition above.

Consider the positive part $v(+)$ and the negative part $v(-)$ of the vector $v$.
The vector $v(+)$ is defined as the componentwise maximum of $v$ and the zero vector,
and the vector $v(-)$ is defined as the componentwise minimum.
Evidently, we have $v=v(+)+v(-)$.
Thus we have two nonnegative vectors $v(+)=u$ and $v(-)=w$
with equal sums of coordinates.

Then classical definition is
$$
\bar \rho(u,w)(=\|u-w\|_{\rho})=\min_{\psi\in \Psi}\sum \rho(x,y)\psi_{x,y},
$$
where $\Psi=\{\psi_{x,y}\}$ is the convex set of all matrices $\psi$ with the following
properties:
$$
\sum_x \psi_{x,y}=u^x;\quad \sum_y \psi_{x,y}=w^y,\quad  \psi_{x,y}\geq 0,\quad  x,y \in X.
$$
Here $u^x$ is the coordinate of the vector $u$ corresponding to the point $x\in X$,
and similarly for $w$.

For the history of the Kantorovich metric, see \cite{V1} and references therein.
For some further development and applications of the finite-dimensional geometry of
this metric, see \cite{V2}.

The last question is also of Arnold's style (see e.g. \cite{A}):
in our definition of the Kantorovich metric, we used
only the root system $A_n$.
My question is: what are the analogs of this definition
(and maybe even of the transportation problem!)
for other series of Lie algebras and other Weyl groups.

\newpage

\begin{figure}[h]
\centering
\includegraphics[height=0.5cm]{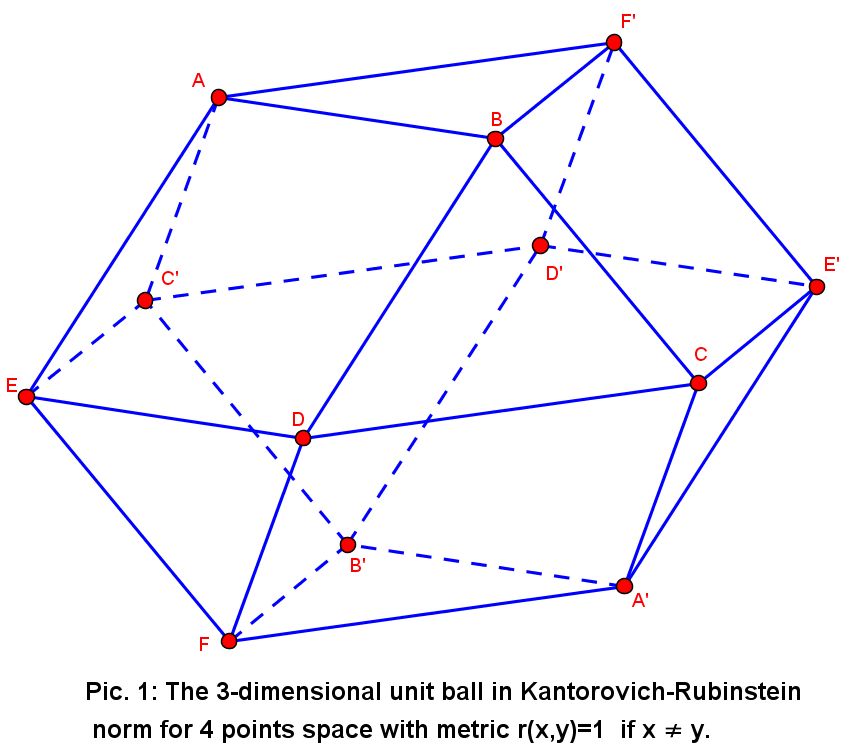}
\label{fig.0}
\end{figure}

\begin{figure}[h]

\bigskip
\bigskip
\bigskip
\bigskip
\bigskip
\bigskip
\bigskip
\bigskip
\bigskip
\bigskip
\bigskip
\bigskip
\bigskip
\bigskip
\bigskip
\bigskip
\bigskip
\bigskip
\bigskip
\bigskip
\bigskip
\bigskip
\bigskip
\bigskip
\bigskip
\bigskip
\bigskip
\bigskip
\bigskip
\bigskip
\bigskip

\centering
\includegraphics[height=1cm]{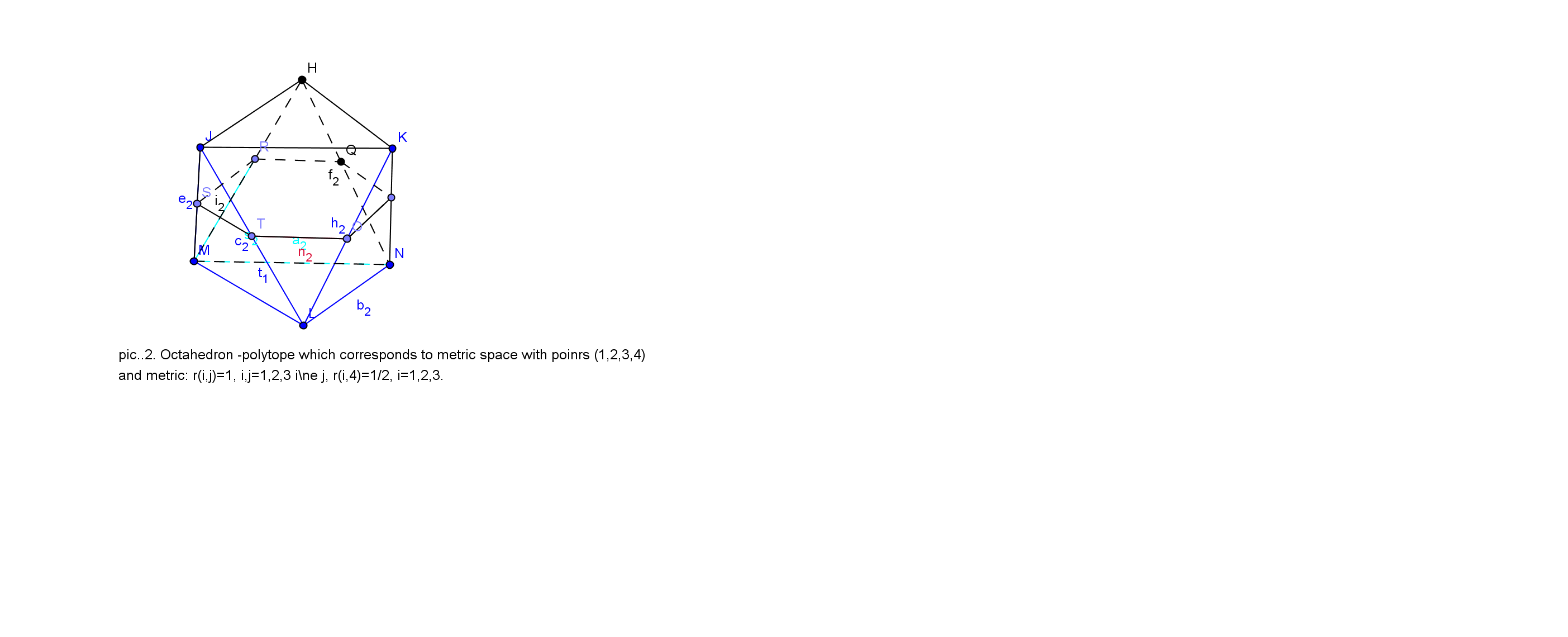}
\label{fig.0}
\end{figure}


\end{document}